\documentclass[11pt,a4paper]{amsart}

\usepackage{amssymb,amsmath,graphics,verbatim}
\usepackage{latexsym}
\usepackage{eucal}
\usepackage{a4wide}

\newtheorem{proposition}{Proposition}[section]
\newtheorem{theorem}{Theorem}[section]
\newtheorem{lemma}[theorem]{Lemma}
\newtheorem{prop}[theorem]{Proposition}
\newtheorem{coro}[theorem]{Corollary}
\newtheorem{remark}[theorem]{Remark}

\newtheorem{definition}[theorem]{Definition}

\newcommand{\mc}{\mathcal}

\newcommand{\rr}{\mathbb{R}}
\newcommand{\R}{\mathbb{R}}
\newcommand{\nn}{\mathbb{N}}
\newcommand{\cc}{\mathbb{C}}

\newcommand{\zz}{\mathbb{Z}}

\newcommand{\eps}{\epsilon}

\newcommand{\pl}{\partial}
\newcommand{\x}{\times}

\newcommand{\til}{\widetilde}
\newcommand{\bbar}{\overline}

\newcommand{\cjd}{\rangle}
\newcommand{\cjg}{\langle}

\newcommand{\demi}{\frac{1}{2}}

\newcommand{\tra}{\textrm{Tr}}

\newcommand{\la}{\lambda}

\def\qed{\hfill$\square$\medskip}

\begin{document}
\title[Inverse scattering]{Inverse scattering at fixed energy on surfaces with Euclidean ends}
\author{Colin Guillarmou}
\address{DMA, U.M.R. 8553 CNRS\\
Ecole Normale Sup\'erieure,\\
45 rue d'Ulm\\ 
F 75230 Paris cedex 05 \\France}
\email{cguillar@dma.ens.fr}

\author{Mikko Salo}
\address{Department of Mathematics \\
University of Helsinki \\
PO Box 68, 00014 Helsinki, Finland}
\email{mikko.salo@helsinki.fi}

\author{Leo Tzou}
\address{Department of Mathematics\\
Stanford University\\
Stanford, CA 94305, USA}
\email{ltzou@math.stanford.edu}

\begin{abstract}
On a fixed Riemann surface  $(M_0,g_0)$ with $N$ Euclidean ends and genus $g$, we show that, under a topological condition,   
the scattering matrix $S_V(\la)$ at frequency $\la > 0$ for the operator $\Delta+V$ 
determines the potential $V$ if $V\in C^{1,\alpha}(M_0)\cap e^{-\gamma d(\cdot,z_0)^j}L^\infty(M_0)$ 
for all $\gamma>0$ and for some $j\in\{1,2\}$, where $d(z,z_0)$ denotes the distance from $z$ 
to a fixed point $z_0\in M_0$. The topological condition is given by $N\geq \max(2g+1,2)$ for $j=1$
and by $N\geq g+1$ if $j=2$. 
In $\rr^2$ this implies that the operator $S_V(\la)$ 
determines any $C^{1,\alpha}$ potential $V$ such that $V(z)=O(e^{-\gamma|z|^2})$ for all $\gamma>0$.
\end{abstract}

\maketitle

\begin{section}{Introduction}
The purpose of this paper is to prove an inverse scattering result at fixed frequency $\la>0$ in dimension $2$. The typical question 
one can ask is to show that the scattering matrix $S_V(\la)$ for the Schr\"odinger operator $\Delta +V$ determines the potential.
This is known to be false if $V$ is only assumed to be Schwartz, by the example of Grinevich-Novikov \cite{GrNo}, 
but  it is also known to be 
true for exponentially decaying potentials (i.e.~$V\in e^{-\gamma|z|}L^\infty(\rr^2)$ for some $\gamma>0$) with norm smaller than a constant depending 
on the frequency $\la$, see Novikov \cite{Nov}. 
For other partial results we refer to \cite{E}, \cite{IS}, \cite{SunU1}, \cite{SunU2}, \cite{SunU3}. The determinacy of $V$ from $S_V(\la)$ when $V$ is compactly supported, without any smallness assumption on the norm, follows from the recent work of Bukhgeim \cite{Bu} on the inverse boundary problem after a standard reduction to the Dirichlet-to-Neumann operator
on a large sphere (see \cite{UhlAsterisque} for this reduction). 

In dimensions $n \geq 3$, it is proved in Novikov \cite{Nov2} (see also \cite{ER} for the case of magnetic Schr\"odinger operators)  
that the scattering matrix at a fixed frequency $\la$ determines an exponentially decaying 
potential. When $V$ is compactly supported this also follows directly from the result by Sylvester-Uhlmann \cite{SylUhl} on the inverse boundary problem, by reducing
 to the Dirichlet-to-Neumann 
operator on a large sphere. 
Melrose \cite{MelStanford} gave a direct proof of the last result based on the methods of \cite{SylUhl}, and this proof was extended to exponentially decaying potentials in \cite{UhlVa} and to the magnetic case in \cite{PSU}.
In the geometric scattering setting, \cite{JSB1,JSB} reconstruct the asymptotic expansion of a potential or metrics
from the scattering operator at fixed frequency on asymptotically Euclidean/hyperbolic manifolds. Further results of this type are given in \cite{W,WY}. 

The method for proving the determinacy of $V$ from $S_V(\la)$ in \cite{MelStanford,UhlVa} is based on the construction of complex geometric optics solutions 
$u(z)=e^{\rho.z}(1+r(\rho,z))$ of $(\Delta+V-\la^2)u=0$ with $\rho\in\cc^n, z\in\rr^n$, and the density of the oscillating scattering solutions 
$u_{\rm sc}(z)=\int_{S^{n-1}} \Phi_V(\la,z,\omega)f(\omega)d\omega$ within those complex geometric optics solutions,  
where $\Phi_V(\la,z,\omega)= e^{i\la \omega.z} + e^{-i\la \omega.z}|z|^{-\demi(n-1)}a(\la,z,\omega)$ are the perturbed plane wave solutions 
(here $\omega\in S^{n-1}$ and $a\in L^\infty$). Unlike when $n \geq 3$, the problem in dimension $2$ is that the set of complex geometrical optics solutions of this type is not  large enough to show that the Fourier transform of $V_1-V_2$ is $0$. 

The real novelty in the recent work of Bukhgeim \cite{Bu} in dimension $2$ is the construction of new complex geometric optics solutions (at least on a bounded domain $\Omega\subset \cc$) of $(\Delta+V_i)u_i=0$ of the form 
$u_1=e^{\Phi/h}(1+r_1(h))$ and $u_2=e^{-\Phi/h}(1+r_2(h))$  with  $0<h\ll 1$ where $\Phi$ is a holomorphic function in $\cc$ with a unique non-degenerate critical point at a fixed 
$z_0\in \cc$ (for instance $\Phi(z)=(z-z_0)^2$), and $||r_j(h)||_{L^p}$ is small as $h\to 0$ for $p>1$. 
These solutions allow to use stationary phase at $z_0$ to get 
\[ \int_{\Omega} (V_1-V_2)u_1\bbar{u_2}  =C (V_1(z_0)-V_2(z_0))h+o(h) , \quad C\not =0\] 
as $h\to 0$ and thus, if the Dirichlet-to-Neumann operators on $\pl\Omega$ are the same, then $V_1(z_0)=V_2(z_0)$. 

One of the problems to extend this to inverse scattering is that a holomorphic function in $\cc$ with a non-degenerate critical point
needs to grow at least quadratically at infinity, which would somehow force to consider potentials $V$ having  Gaussian decay. On the other hand,
if we allow the function to be meromorphic with simple poles, then we can construct such functions, having a single critical point at any given point 
$p$, for instance by considering $\Phi(z)=(z-p)^2/z$. Of course, with such $\Phi$ we then need to work on 
$\cc\setminus\{0\}$, which is conformal to a surface with no hole but with $2$ Euclidean ends, and $\Phi$ has linear growth in the 
ends. In general, on a surface with genus $g$ and $N$ Euclidean ends, we can use the Riemann-Roch theorem to construct holomorphic 
functions with linear or quadratic growth in the ends, the dimension of the space of such functions depending  on $g,N$. 

In the present work, we apply this idea to obtain an inverse scattering result for $\Delta_{g_0}+V$ on a fixed Riemann surface $(M_0,g_0)$
with Euclidean ends, under some topological condition on $M_0$ and some decay condition on $V$.
\begin{theorem}\label{mainth}
Let $(M_0,g_0)$ be a non-compact Riemann surface with genus $g$ and $N$ ends isometric to $\rr^2\setminus\{|z|\leq 1\}$ 
with metric $|dz|^2$. Let $V_1$ and $V_2$ be two potentials in $C^{1,\alpha}(M_0)$ with $\alpha>0$, and such that 
 $S_{V_1}(\la)=S_{V_2}(\la)$ for some  $\la > 0$. Let $d(z,z_0)$ denote the distance  
between $z$ and a fixed point $z_0\in M_0$.\\ 
(i) If  $N\geq \max(2g+1,2)$ and $V_i\in e^{-\gamma d(\cdot,z_0)}L^\infty(M_0)$ for all $\gamma>0$, then $V_1=V_2$.\\
(ii) If $N\geq g+1$ and $V_i\in e^{-\gamma d(\cdot,z_0)^2}L^\infty(M_0)$ for all $\gamma>0$, then $V_1=V_2$.
\end{theorem}

In $\rr^2$, where $g = 0$ and $N = 1$, we have an immediate corollary: 
\begin{coro}\label{th2}
Let $\la > 0$ and let $V_1,V_2\in C^{1,\alpha}(\rr^2) \cap e^{-\gamma |z|^2}L^\infty(\rr^2)$ for all $\gamma>0$. If the scattering matrices satisfy $S_{V_1}(\la)=S_{V_2}(\la)$, then $V_1=V_2$.
\end{coro}
This is an improvement on the result of Bukhgeim
\cite{Bu} which shows identifiability for compactly supported functions, and in a certain sense on the result of Novikov \cite{Nov}
since it is assumed there that the potential has to be of small $L^\infty$ norm. 

The structure of the paper is as follows. In Section \ref{sec_holomorphic} we employ the Riemann-Roch theorem and a transversality argument to construct Morse holomorphic functions on $(M_0,g_0)$ with linear or quadratic growth in the ends. Section \ref{Carleman} considers Carleman estimates with harmonic weights on $(M_0,g_0)$, where suitable convexification and weights at the ends are required since the surface is non compact. Complex geometrical optics solutions are constructed in Section \ref{CGOriemann}.
Section \ref{sec_scattering} discusses direct scattering theory on surfaces with Euclidean ends and contains the proof that scattering solutions are dense in the set of suitable solutions, and Section \ref{sec_identify} gives the proof of Theorem \ref{mainth}. Finally, there is an appendix discussing a Paley-Wiener type result for functions with Gaussian decay which is needed to prove density of scattering solutions.

\subsection*{Acknowledgements}
M.S.~is supported partly by the Academy of Finland. C.G. is supported by ANR grant ANR-09-JCJC-0099-01, and is grateful to the Math. Dept. of Helsinki where part of this work was done. 

\end{section}

\begin{section}{Holomorphic Morse functions on a surface with Euclidean ends} \label{sec_holomorphic}

\subsection{Riemann surfaces with Euclidean ends}
Let $(M_0,g_0)$ be a non-compact connected smooth Riemannian surface with $N$ ends $E_1,\dots,E_N$ which are Euclidean, i.e. isometric to 
$\cc\setminus \{|z|\le 1\}$ with metric $|dz|^2$. By using a complex inversion $z\to 1/z$, each end is also isometric to a pointed disk
\[ E_i \simeq \{|z|\leq 1, z\not=0\} \textrm{ with metric }\frac{|dz|^2}{|z|^4}\]
thus conformal to the Euclidean metric on the pointed disk. The surface $M_0$ can then be compactified by adding the points corresponding  
to $z=0$ in each pointed disk corresponding to an end $E_i$, we obtain a closed Riemann surface $M$ with a natural complex structure induced by that of 
$M_0$, or equivalently a smooth conformal class on $M$ induced by that of $M_0$. Another way of thinking is to say that $M_0$ is the closed Riemann surface $M$ with $N$ points $e_1,\dots,e_N$ removed.   
The Riemann surface $M$ has  holomorphic charts $z_\alpha:U_{\alpha}\to \cc$ and we will denote by $z_1,\dots z_N$ the complex coordinates 
corresponding to the ends of $M_0$, or equivalently to the neighbourhoods of the points $e_i$.
 The Hodge star operator $\star$ acts on the cotangent bundle $T^*M$, its eigenvalues are 
$\pm i$ and the respective eigenspaces $T_{1,0}^*M:=\ker (\star+i{\rm Id})$ and $T_{0,1}^*M:=\ker(\star -i{\rm Id})$
are sub-bundles of the complexified cotangent bundle $\cc T^*M$ and the splitting $\cc T^*M=T^*_{1,0}M\oplus T_{0,1}^*M$ holds as complex vector spaces.
Since $\star$ is conformally invariant on $1$-forms on $M$, the complex structure depends only on the conformal class of $g$.
In holomorphic coordinates $z=x+iy$ in a chart $U_\alpha$,
one has $\star(udx+vdy)=-vdx+udy$ and
\[T_{1,0}^*M|_{U_\alpha}\simeq \cc dz ,\quad T_{0,1}^*M|_{U_\alpha}\simeq \cc d\bar{z}  \]
where $dz=dx+idy$ and $d\bar{z}=dx-idy$. We define the natural projections induced by the splitting of $\cc T^*M$ 
\[\pi_{1,0}:\cc T^*M\to T_{1,0}^*M ,\quad \pi_{0,1}: \cc T^*M\to T_{0,1}^*M.\]
The exterior derivative $d$ defines the de Rham complex $0\to \Lambda^0\to\Lambda^1\to \Lambda^2\to 0$ where $\Lambda^k:=\Lambda^kT^*M$
denotes the real bundle of $k$-forms on $M$. Let us denote $\cc\Lambda^k$ the complexification of $\Lambda^k$, then
the $\pl$ and $\bar{\pl}$ operators can be defined as differential operators 
$\pl: \cc\Lambda^0\to T^*_{1,0}M$ and $\bar{\pl}:\cc\Lambda0\to T_{0,1}^*M$ by 
\[\pl f:= \pi_{1,0}df ,\quad \bar{\pl}f:=\pi_{0,1}df,\]
they satisfy $d=\pl+\bar{\pl}$ and are expressed in holomorphic coordinates by
\[\pl f=\pl_zf\, dz ,\quad \bar{\pl}f=\pl_{\bar{z}}f \, d\bar{z},\]  
with $\pl_z:=\demi(\pl_x-i\pl_y)$ and $\pl_{\bar{z}}:=\demi(\pl_x+i\pl_y)$.
Similarly, one can define the $\pl$ and $\bar{\pl}$ operators from $\cc \Lambda^1$ to $\cc \Lambda^2$ by setting 
\[\pl (\omega_{1,0}+\omega_{0,1}):= d\omega_{0,1}, \quad \bar{\pl}(\omega_{1,0}+\omega_{0,1}):=d\omega_{1,0}\]
if $\omega_{0,1}\in T_{0,1}^*M$ and $\omega_{1,0}\in T_{1,0}^*M$.
In coordinates this is simply
\[\pl(udz+vd\bar{z})=\pl v\wedge d\bar{z},\quad \bar{\pl}(udz+vd\bar{z})=\bar{\pl}u\wedge d{z}.\]
If $g$ is a metric on $M$ whose conformal class induces the complex structure $T_{1,0}^*M$, 
there is a natural operator, the Laplacian acting on functions and defined by 
\[\Delta f:= -2i\star \bar{\pl}\pl f =d^*d \]
where $d^*$ is the adjoint of $d$ through the metric $g$ and $\star$ is the Hodge star operator mapping 
$\Lambda^2$ to $\Lambda^0$ and induced by $g$ as well.
 
\subsection{Holomorphic functions}
We are going to construct Carleman weights given by holomorphic functions on $M_0$ which grow at most linearly or quadratically 
in the ends.  We will use the Riemann-Roch theorem, following ideas of \cite{GT}, however, the difference in the present case
is that we have very little freedom to construct these holomorphic functions, simply because there is just a finite dimensional space of such 
functions by Riemann-Roch.  
For the convenience of the reader, and to fix notations, we recall the usual Riemann-Roch index theorem 
(see Farkas-Kra \cite{FK} for more details). A divisor $D$ on $M$ is an element 
\[D=\big((p_1,n_1), \dots, (p_k,n_k)\big)\in (M\x\zz)^k, \textrm{ where }k\in\nn\]  
which will also be denoted $D=\prod_{i=1}^kp_i^{n_i}$ or $D=\prod_{p\in M}p^{\alpha(p)}$ where $\alpha(p)=0$
for all $p$ except $\alpha(p_i)=n_i$. The inverse divisor of $D$ is defined to be 
$D^{-1}:=\prod_{p\in M}p^{-\alpha(p)}$ and 
the degree of the divisor $D$ is defined by $\deg(D):=\sum_{i=1}^kn_i=\sum_{p\in M}\alpha(p)$. 
A non-zero meromorphic function on $M$ is said to have divisor $D$ if $(f):=\prod_{p\in M}p^{{\rm ord}(p)}$ is equal to $D$,
where ${\rm ord}(p)$ denotes the order of $p$ as a pole or zero of $f$ (with positive sign convention for zeros). Notice that 
in this case we have $\deg(f)=0$.
For divisors $D'=\prod_{p\in M}p^{\alpha'(p)}$ and $D=\prod_{p\in M}p^{\alpha(p)}$, 
we say that $D'\geq D$ if $\alpha'(p)\geq \alpha(p)$ for all $p\in M$.
The same exact notions apply for meromorphic $1$-forms on $M$. Then we define for a divisor $D$
\[\begin{gathered}
r(D):=\dim (\{f \textrm{ meromorphic function on } M; (f)\geq D\}\cup\{0\}),\\
i(D):=\dim(\{u\textrm{ meromorphic 1 form on }M; (u)\geq D\}\cup\{0\}).
\end{gathered}\]
The Riemann-Roch theorem states the following identity: for any divisor $D$ on the closed Riemann surface $M$ of genus $g$, 
\begin{equation}\label{riemannroch}
r(D^{-1})=i(D)+\deg(D)-g+1.
\end{equation}
Notice also that for any divisor $D$ with $\deg(D)>0$, one has $r(D)=0$ since $\deg(f)=0$ for all $f$ meromorphic. 
By \cite[Th. p70]{FK}, let $D$ be a divisor, then for any non-zero meromorphic 1-form $\omega$ on $M$, one has 
\begin{equation}\label{abelian}
i(D)=r(D(\omega)^{-1})
\end{equation}
which is thus independent of $\omega$. 
For instance, if $D=1$, we know that the only holomorphic function on $M$ is $1$ and 
one has $1=r(1)=r((\omega)^{-1})-g+1$ and thus $r((\omega)^{-1})=g$ if $\omega$ is a non-zero meromorphic $1$ form. Now 
if $D=(\omega)$, we obtain again from \eqref{riemannroch}
\[ g=r((\omega)^{-1})=2-g+\deg((\omega))\] 
 which gives $\deg((\omega))=2(g-1)$ for any non-zero meromorphic $1$-form $\omega$. In particular, if $D$ is a divisor such that
$\deg(D)>2(g-1)$, then  we get 
$\deg(D(\omega)^{-1})=\deg(D)-2(g-1)>0$ and thus $i(D)=r(D(\omega)^{-1})=0$, which implies by \eqref{riemannroch}
\begin{equation}\label{divisors}
\deg(D)>2(g-1)\Longrightarrow r(D^{-1})= \deg(D)-g+1\geq g.
\end{equation}
Now we deduce the 
\begin{lemma}\label{holofcts} 
Let $e_1,\dots,e_{N}$ be distinct points on a closed Riemann surface $M$ with genus $g$, and let $z_0$ be another 
point of $M\setminus\{e_1,\dots,e_{N}\}$. If $N\geq \max(2g+1,2)$, the following hold true:\\ 
(i) there exists a meromorphic function $f$ on $M$ with at most simple 
poles, all contained in $\{e_1,\dots,e_{N}\}$, such that $\pl f(z_0)\not=0$,\\
(ii) there exists a meromorphic function $h$
 on $M$ with at most simple poles, all contained in $\{e_1,\dots,e_{N}\}$, 
such that $z_0$ is a zero of order at least $2$ of $h$.
\end{lemma}
\noindent\textsl{Proof}. Let first $g \geq 1$, so that $N \geq 2g+1$. By the discussion before the Lemma, we know that there are at least $g+2$ linearly independent (over $\cc$) 
meromorphic functions $f_0,\dots,f_{g+1}$ on $M$ with at most simple poles, all contained in $\{e_1,\dots,e_{2g+1}\}$. Without loss of generality,
one can set $f_0=1$ and by linear combinations we can assume that $f_1(z_0)=\dots=f_{g+1}(z_0)=0$. Now consider the divisor  
$D_j=e_1\dots e_{2g+1}z_0^{-j}$ for $j=1,2$, with degree $\deg(D_j)=2g+1-j$, then by the Riemann-Roch formula (more precisely \eqref{divisors})
\[ r(D_j^{-1})=g+2-j.\]
Thus, since $r(D_1^{-1})>r(D_2^{-1})=g$ and using the assumption that $g \geq 1$, we deduce that there is a function in ${\rm span}(f_1,\dots,f_{g+1})$ 
which has a zero of order $2$ at $z_0$ and a function which has a zero of order exactly $1$ at $z_0$. The same method clearly 
works if $g=0$ by taking two points $e_1,e_2$ instead of just $e_1$. 
\qed
 
If we allow double poles instead of simple poles, the proof of Lemma \ref{holofcts} shows the
\begin{lemma}\label{doublepoles}
Let $e_1,\dots,e_{N}$ be distinct points on a closed Riemann surface $M$ with genus $g$, and let $z_0$ be another 
point of $M\setminus\{e_1,\dots,e_{N}\}$. If $N\geq g+1$, then the following hold true:\\ 
(i) there exists a meromorphic function $f$ on $M$ with at most double 
poles, all contained in $\{e_1,\dots,e_{N}\}$, such that $\pl f(z_0)\not=0$,\\
(ii) there exists a meromorphic function $h$
 on $M$ with at most double poles, all contained in $\{e_1,\dots,e_{N}\}$, 
such that $z_0$ is a zero of order at least $2$ of $h$.\\
\end{lemma}  
 
\subsection{Morse holomorphic functions with prescribed critical points} \label{morseholo}
We follow in this section the arguments used in \cite{GT} to construct holomorphic functions with non-degenerate critical points (i.e.~Morse holomorphic functions) on 
the surface $M_0$ with genus $g$ and $N$ ends, such that these functions have at most linear growth (resp.~quadratic growth) 
in the ends if $N\geq \max(2g+1,2)$ (resp.~if $N\geq g+1$).
We let $\mc{H}$ be the complex  vector space spanned by the meromorphic functions on $M$ with divisors larger or equal to  
$e_1^{-1}\dots e_{N}^{-1}$ (resp. by $e_1^{-2}\dots e_{N}^{-2}$) if we work with functions having linear growth (resp.~quadratic growth), 
where $e_1,\dots e_{N}\in M$ are points corresponding to the ends of $M_0$ as explained in
Section \ref{sec_holomorphic}.  
Note that $\mc{H}$ is a complex vector space of complex dimension greater or equal to $N-g+1$ (resp.~$2N-g+1$) 
for the $e_1^{-1}\dots e_{N}^{-1}$ divisor (resp.~the $e_1^{-2}\dots e_{N}^{-2}$ divisor).
We will also consider the real vector space $H$ spanned by the real parts and imaginary parts of functions in $\mc{H}$, this is a real vector space 
which  admits a  Lebesgue measure. We now prove the following 
\begin{lemma}\label{morsedense}
The set of functions $u\in H$ which are not Morse in $M_0$ has measure $0$ in $H$, in particular its complement is dense in $H$.
\end{lemma}
\noindent{\bf Proof}. We use an argument very similar to that used by Uhlenbeck \cite{Uh}.
We start by defining $m: M_0\times H\to T^*M_0$ by $(p,u) \mapsto (p,du(p))\in T_p^*M_0$. 
This is clearly a smooth map, linear in the second variable, moreover $m_u:=m(.,u)=(\cdot, du(\cdot))$ is  
smooth on $M_0$. The map $u$ is a Morse function if and only if 
$m_u$ is transverse to the zero section, denoted $T_0^*M_0$, of $T^*M_0$, i.e.~if 
\[\textrm{Image}(D_{p}m_u)+T_{m_u(p)}(T_0^*M_0)=T_{m_u(p)}(T^*M_0),\quad \forall p\in M_0 \textrm{ such that }m_u(p)=(p,0).\]
This is equivalent to the fact that the Hessian of $u$ at critical points is 
non-degenerate (see for instance Lemma 2.8 of \cite{Uh}). 
We recall the following transversality result, the proof of which is contained in \cite[Th.2]{Uh} by replacing Sard-Smale theorem by the usual 
finite dimensional Sard theorem:
\begin{theorem}\label{transv}
Let $m : X\times H \to W$ be a $C^k$ map and $X, W$ be smooth manifolds and $H$  a finite dimensional vector space, 
if  $W'\subset W$ is a submanifold such that $k>\max(1,\dim X-\dim W+\dim W')$, then the transversality of the map $m$ to $W'$  implies that 
the complement of the set $\{u\in H; m_u \textrm{ is transverse to } W'\}$ in $H$ has Lebesgue measure $0$.
\end{theorem} 
We want to apply this result with $X:=M_0$, $W:=T^*M_0$ and $W':=T^*_0M_0$, and with the map $m$ as defined above. 
We have thus proved our Lemma if one can show that $m$ is transverse to $W'$. 
Let $(p,u)$ such that $m(p,u)=(p,0)\in W'$. Then identifying $T_{(p,0)}(T^*M_0)$ with $T_pM_0\oplus T^*_pM_0$, one has
\[Dm_{(p,u)}(z,v)=(z,dv(p)+{\rm Hess}_p(u)z)\]
where ${\rm Hess}_p(u)$ is the Hessian of $u$ at the point $p$, viewed as a linear map from $T_pM_0$ to $T^*_pM_0$ (note that this is different from the covariant Hessian defined by the Levi-Civita connection). 
To prove that $m$ is transverse to $W'$ we need to show that $(z,v)\to (z, dv(p)+{\rm Hess}_p(u)z)$ is onto from $T_pM_0\oplus H$ 
to $T_pM_0\oplus T^*_pM_0$, which is realized if the map $v\to dv(p)$ from $H$ to  $T_p^*M_0$ is onto.
But from Lemma \ref{holofcts}, we know that there exists a meromorphic function $f$ with real part $v={\rm Re}(f)\in H$ 
such that $v(p)=0$ and $dv(p)\not=0$ as an element of $T^*_pM_0$. We can then take $v_1:=v$ and $v_2:={\rm Im}(f)$, which are 
functions of $H$ such that $dv_1(p)$ and $dv_2(p)$ are
linearly independent in $T^*_pM_0$ by the Cauchy-Riemann equation $\bar{\pl} f=0$. 
This shows our claim and ends the proof by using Theorem \ref{transv}.\qed

In particular, by the Cauchy-Riemann equation, this Lemma implies that the set of Morse 
functions in $\mc{H}$ is dense in $\mc{H}$.  We deduce 
\begin{proposition}
\label{criticalpoints}
There exists a dense set of points $p$ in $M_0$ such that there exists a Morse holomorphic function $f\in\mc{H}$ on $M_0$ 
which has a critical point at $p$.
\end{proposition} 
\noindent\textsl{Proof}.  
Let $p$ be a point of $M_0$ and let $u$ be a holomorphic function with a zero of order at least $2$ at $p$, 
the existence is ensured by Lemma \ref{holofcts}. 
Let $B(p,\eta)$ be a any small ball of radius $\eta>0$ near $p$, then by 
Lemma \ref{morsedense}, for any $\eps>0$, we can approach $u$ by a holomorphic Morse function $u_\eps\in\mc{H}_\eps$ 
which is at distance less than $\eps$ of $u$ in a fixed norm on the finite dimensional space $\mc{H}$. 
Rouch\'e's theorem for $\pl_z u_\eps$ and $\pl_zu$ (which are viewed as functions locally near $p$) 
implies that $\pl_z u_{\eps}$ has at least one zero of order exactly $1$ in $B(p,\eta)$ if $\eps$ is chosen small enough. 
Thus there is a Morse function in $\mc{H}$ with a critical point  arbitrarily close to $p$.\qed
\end{section}

\begin{remark}
In the case where the surface $M$ has genus $0$ and $N$ ends,  we have an explicit formula for the function in Proposition \ref{criticalpoints}: indeed 
$M_0$ is conformal to $\cc\setminus \{e_1,\dots, e_{N-1}\}$ for some
$e_i\in\cc$ - i.e.~the Riemann sphere minus $N$ points - then the function $f(z)=(z-z_0)^2/(z-e_1)$ with $z_0\not\in\{e_1,\dots,e_{N-1}\}$ 
has $z_0$ for unique critical point in 
$\cc\setminus \{e_1,\dots, e_{N-1}\}$ and it is non-degenerate. 
\end{remark}

We end this section by the following Lemmas which will be used for the amplitude of the complex geometric optics solutions but not for the phase.
\begin{lemma}\label{amplitude}
For any $p_0, p_1,\dots p_n\in M_0$ some points of $M_0$ and $L\in\nn$, then there exists a function $a(z)$ holomorphic on $M_0$ 
which vanishes to order $L$ at all $p_j$ for $j=1,\dots,n$ and such that $a(p_0)\not=0$.
 Moreover $a(z)$ can be chosen to have  at most polynomial growth in the ends, i.e.
$|a(z)|\leq C|z|^{J}$ for some $J\in\nn$.
\end{lemma}
\noindent\textsl{Proof}. It suffices to find on $M$ some meromorphic function with divisor greater or equal to  
$D:=e_1^{-J}\dots e_N^{-J}p_1^L\dots p_n^{L}$ but not greater or equal to $Dp_0$ 
and this is insured by Riemann-Roch theorem as long as  $JN-nL\geq 2g$ since then 
$r(D)=-g+1+JN-nL$ and $r(Dp_0)=-g+JN-nL$.
\qed

\begin{lemma}
\label{control the zero}
Let $\{p_0, p_1,..,p_n\}\subset M_0$ be a set of $n+1$ disjoint points. Let $c_0,c_1,\dots, c_K\in\cc$, $L\in\nn$, 
and let $z$ be a complex coordinate near $p_0$ such that $p_0=\{z=0\}$. 
Then  there exists a holomorphic function $f$ on $M_0$ with zeros of order at least $L$ at each $p_j$,
such that $f(z)=c_0 + c_1z +...+ c_K z^K+ O(|z|^{K+1})$ in the coordinate $z$.  Moreover $f$ 
can be chosen so that there is $J\in\nn$ such that, in the ends, $|\pl_z^\ell f(z)|=O (|z|^{J})$ for all $\ell\in\nn_0$.
\end{lemma}
\noindent\textsl{Proof}. The proof goes along the same lines as in Lemma \ref{amplitude}. By induction on $K$ and linear combinations, 
it suffices to prove it for $c_0=\dots=c_{K-1}=0$. As in the proof of Lemma \ref{amplitude}, if $J$ is taken large enough, there exists 
a function with divisor greater or equal to $D:=e_1^{-J}\dots e_N^{-J}p_0^{K-1}p_1^L\dots p_n^{L}$ but not greater or equal 
to $Dp_0$. Then it suffices to multiply this function by $c_K$ times the inverse of 
the coefficient of $z^K$ in its Taylor expansion at $z=0$.
\qed

\subsection{Laplacian on weighted spaces} 
Let $x$ be a smooth positive function on $M_0$, which is equal to $|z|^{-1}$ for $|z|>r_0$ in the ends 
$E_i\simeq \{z\in\cc; |z|>1\}$, where $r_0$ is a large fixed number. 
We now show that the Laplacian $\Delta_{g_0}$ on a surface with Euclidean ends has a right inverse on the weighted spaces 
$x^{-J}L^2(M_0)$ for $J\notin\nn$ positive.
\begin{lemma}\label{rightinv}
For any $J>-1$  which is not an integer, there exists a continuous operator $G$ mapping $x^{-J}L^2(M_0)$ to $x^{-J-2}L^2(M_0)$  
such that $\Delta_{g_0}G={\rm Id}$.
\end{lemma}
\noindent\textsl{Proof}. Let $g_b:=x^2g_0$ be a metric conformal to $g_0$. The metric $g_b$ in the ends can be written 
$g_b=dx^2/x^2+d\theta_{S^1}^2$ by using radial coordinates $x=|z|^{-1},\theta=z/|z|\in S^1$,  
this is thus a b-metric in the sense of Melrose \cite{APS}, giving the surface a geometry of surface with cylindrical ends. Let us define for $m\in\nn_0$
\[H^m_b(M_0):=\{u\in L^2(M_0;{\rm dvol}_{g_b}); (x\pl_x)^j \pl_\theta^ku\in L^2(M_0;{\rm dvol}_{g_b}) \textrm{ for all }j+k\leq m\}.\] 
The Laplacian has the form $\Delta_{g_b}=-(x\pl_x)^2+\Delta_{S^1}$ in the ends, and 
the indicial roots of $\Delta_{g_b}$ in the sense of Section 5.2 of \cite{APS} are given by the complex numbers $\la$ such that
$x^{-i\la}\Delta_{g_b}x^{i\la}$ is not invertible as an operator acting on the circle $S_{\theta}^1$. Thus the indicial roots are the solutions of $\la^2+k^2=0$ 
where $k^2$ runs over the eigenvalues of $\Delta_{S^1}$, that is, $k \in \zz$. The roots are simple at $\pm i k \in i\zz\setminus \{0\}$ and $0$ is a 
double root. In Theorem 5.60 of \cite{APS}, Melrose proves that $\Delta_{g_b}$ is Fredholm on $x^{a}H^2_b(M_0)$ if and only if $-a$ is not the imaginary part of some indicial root, that is here $a\not\in \zz$.  
For $J>0$, the kernel of $\Delta_{g_b}$ on the space $x^{J}H^2_b(M_0)$ is clearly trivial by an energy estimate. Thus $\Delta_{g_b}: x^{-J} H^0_b(M_0) \to x^{-J} H^{-2}_b(M_0)$ is surjective for $J > 0$ and $J \not\in \zz$, and the same then holds for $\Delta_{g_b}: x^{-J} H^2_b(M_0) \to x^{-J} H^{0}_b(M_0)$ by elliptic regularity.

Now we can use Proposition 5.64 of \cite{APS}, which asserts, for all positive $J \not\in \zz$, the existence of a  pseudodifferential operator $G_b$ mapping continuously $x^{-J}H^0_b(M_0)$ to $x^{-J}H^2_b(M_0)$ such  that $\Delta_{g_b}G_{b}={\rm Id}$.
Thus if we set $G=G_bx^{-2}$, we have $\Delta_{g_0}G={\rm Id}$ and $G$ maps continuously 
$x^{-J+1}L^2(M_0)$ to $x^{-J-1}L^2(M_0)$ (note that $L^2(M_0)=xH^0_b(M_0)$). 
\qed

\begin{section}{Carleman Estimate for Harmonic Weights with Critical Points}\label{Carleman}

\subsection{The linear weight case}
In this section, we prove a Carleman estimate using harmonic weights with non-degenerate critical points, in a
way similar to \cite{GT}. Here however we need to work on a non compact surface and with weighted spaces.
We first consider a Morse holomorphic function $\Phi\in \mc{H}$ obtained from Proposition \ref{criticalpoints}
with the condition that $\Phi$ has linear growth in the ends, which corresponds to the case where $V\in e^{-\gamma/x}L^\infty(M_0)$ for all $\gamma>0$. The Carleman weight will be the harmonic function $\varphi := \mathrm{Re}(\Phi)$.
We let  $x$ be a positive smooth function  on $M_0$ such that  $x=|z|^{-1}$  in 
the complex charts $\{z\in\cc; |z|>1\}\simeq E_i$ covering the end $E_i$.\\ 

Let $\delta\in(0,1)$ be small and let us take $\varphi_0\in x^{-\alpha}L^2(M_0)$ a solution of $\Delta_{g_0}\varphi_0=x^{2-\delta}$, a solution exists 
by Proposition \ref{rightinv} if $\alpha>1+\delta$. 
Actually, by using Proposition 5.61 of \cite{APS}, if we choose $\alpha<2$, then it is easy to see that 
$\varphi_0$ is smooth on $M_0$ and has polyhomogeneous expansion as $|z|\to \infty$, with leading asymptotic in the end $E_i$ given by
$\varphi_0= -x^{-\delta}/\delta^2+
c_i\log(x)+ d_i+O(x)$ for some $c_i,d_i$ which are smooth functions in $S^1$.
For $\eps>0$ small, we define the convexified weight $\varphi_\eps:=\varphi-\frac{h}{\eps}\varphi_0$.\\

We recall from the proof of Proposition 3.1 in \cite{GT} the following estimate which is valid in any compact set $K\subset M_0$: for all
$w\in C_0^\infty(K)$, we have 
\begin{equation}\label{carlemaninK}
\frac{C}{\eps}\Big(\frac{1}{h}\|w\|_{L^2}^2 + \frac{1}{h^2}\|w |d\varphi|\|_{L^2}^2 +\frac{1}{h^2}\|w |d\varphi_{\eps}|\|_{L^2}^2 +
\|dw\|_{L^2(K)}^2\Big) \leq \|e^{\varphi_{\eps}/h}\Delta_g e^{-\varphi_{\eps}/h} w\|_{L^2}^2
\end{equation} 
where $C$ depends on $K$ but not on $h$ and $\eps$.\\

So for functions supported in the end $E_i$, it clearly suffices to obtain a Carleman estimate in 
$E_i\simeq \rr^2\setminus\{|z|\leq 1\}$  by using the Euclidean coordinate $z$ of the end. 
\begin{prop}\label{carlemaninend}
Let $\delta\in(0,1)$, and $\varphi_\eps$ as above, then there exists $C>0$ such that for all $\eps\gg h>0$ small enough, and all $u\in C_0^\infty(E_i)$
\[ h^2||e^{\varphi_{\eps}/h}(\Delta-\la^2)e^{-\varphi_{\eps}/h}u||^2_{L^2}\geq \frac{C}{\eps} (||x^{1-\frac{\delta}{2}}u||^2_{L^2}+ 
h^2||x^{1-\frac{\delta}{2}}du||^2_{L^2}).\]
\end{prop}
\noindent\textsl{Proof}.
The metric $g_0$ can be extended to $\rr^2$ to be the Euclidean metric
and we shall denote by $\Delta$ the flat positive Laplacian on $\rr^2$.
Let us write $P:=\Delta_{g_0}-\la^2$, then the operator $P_h:= h^2e^{\varphi_\eps/h}Pe^{-\varphi_\eps/h}$
is given by 
\[ P_h=h^2\Delta-|d\varphi_\eps|^2+2h\nabla\varphi_\eps.\nabla- h\Delta \varphi_\eps-h^2\la^2, \]
following the notation of \cite[Chap. 4.3]{EvZw}, it is a semiclassical operator in  $S^0(\cjg\xi\cjd^2)$ with  semiclassical full Weyl symbol 
\[\sigma(P_h):=|\xi|^2-|d\varphi_\eps|^2-h^2\la^2+2i\cjg d\varphi_\eps,\xi\cjd =a+ib.\]
We can define $A:=(P_h+P_h^*)/2=h^2\Delta-|d\varphi_\eps|^2-h^2\la^2$ and
 $B:=(P_h-P_h^*)/2i=-2ih\nabla\varphi_\eps.\nabla+ih\Delta\varphi_\eps$ which have respective semiclassical full symbols
$a$ and $b$, i.e.  $A={\rm Op}_h(a)$  and $B={\rm Op}_h(b)$ for the Weyl quantization. Notice that $A,B$ are symmetric operators, 
thus for all $u\in C_0^\infty(E_i)$
\begin{equation}\label{A+iB} 
||(A+iB)u||^2=\cjg (A^2+B^2+i[A,B])u,u\cjd .
\end{equation}
It is easy to check that the operator $ih^{-1}[A,B]$ is a semiclassical differential operator in $S^0(\cjg\xi\cjd^2)$ 
with full semiclassical symbol  
\begin{equation}\label{symbol[A,B]}
\{a,b\}(\xi)=  4(D^2\varphi_\eps(d\varphi_\eps,d\varphi_\eps)+D^2\varphi_\eps(\xi,\xi))
\end{equation}

Let us now decompose the Hessian of $\varphi_\eps$ in the basis $(d\varphi_\eps, \theta)$ where $\theta$ is a 
covector orthogonal to $d\varphi_\eps$ and of norm $|d\varphi_\eps|$. This yields coordinates $\xi=\xi_0d\varphi_\eps+\xi_1\theta$ 
and there exist smooth functions $M,N,K$ so that
\[D^2\varphi_\eps(\xi,\xi)=|d\varphi_\eps|^2(M\xi_0^2+N\xi_1^2+2K\xi_0\xi_1).\]
Notice that $\varphi_\eps$ has a polyhomogeneous expansion at infinity of the form 
\[\varphi_\eps(z)= \gamma.z+ \frac{h}{\eps} \frac{r^{\delta}}{\delta^2} +c_1\log(r)+ c_2+c_3 r^{-1}+O(r^{-2}) \]
where $r=|z|,\omega=z/r, \gamma=(\gamma_1,\gamma_2)\in \rr^2$ and $c_i$ are some smooth functions  on $S^1$ depending on $h$; 
in particular we have 
\[ d\varphi_\eps = \gamma_1dz_1+\gamma_2dz_2+O(r^{-1+\delta}), \quad 
 \pl_{z}^\alpha\pl_{\bar{z}}^\beta \varphi_\eps(z)=O(r^{-2+\delta}) \,\,\textrm{ for all }  \alpha+\beta\geq 2\] 
which implies that $M,N,K\in r^{-2+\delta}L^\infty(E_i)$. 
Then one can write 
\[\begin{split} 
\{a,b\}= & 4|d\varphi_\eps|^2(M+M\xi_0^2+N\xi_1^2+2K\xi_0\xi_1)\\
=&4(N(a+h^2\la^2) +((M-N)\xi_0+2K\xi_1)b/2+(N+M)|d\varphi_\eps|^2)
\end{split}\]
and since $M+N=\tra(D^2\varphi_\eps)=-\Delta\varphi_\eps=h\Delta\varphi_0/\eps$ we obtain 
\begin{equation}\label{{a,b}}\begin{gathered}
\{a,b\}=4|d\varphi_\eps|^2(c(z)(a+h^2\la^2) +\ell(z,\xi)b+\frac{h}{\eps} r^{-2+\delta}),\\
 c(z)=\frac{N}{|d\varphi_\eps|^2}, \,\,\, \ell(z,\xi)=\frac{(M-N)\xi_0+2K\xi_1}{2|d\varphi_\eps|^2}.
\end{gathered}
\end{equation} 
Now, we take a smooth extension of $|d\varphi_\eps|^2, a(z,\xi),\ell(z,\xi)$ and $r$ to $z\in\rr^2$, this can done for instance by  extending 
$r$ as a smooth positive function on $\rr^2$ and  then extending $d\varphi$ and $d\varphi_0$ to smooth non vanishing $1$-forms 
on $\rr^2$ (not necessarily exact) so that $|d\varphi_\eps|^2$ is smooth positive (for small $h$) and polynomial in $h$
and $a,\ell$ are of the same form as in $\{|z|>1\}$.
Let us define the symbol and quantized differential operator on $\rr^2$
\[ e:=4|d\varphi_\eps|^2(c(z)(a+h^2\la^2) +\ell(z,\xi)b), \quad E:={\rm Op}_h(e)\] 
and write 
\begin{equation}\label{simpl}
\begin{gathered}
ih^{-1}r^{1-\frac{\delta}{2}}[A,B]r^{1-\frac{\delta}{2}}= hF + r^{1-\frac{\delta}{2}}Er^{1-\frac{\delta}{2}}-\frac{h}{\eps}(A^2+B^2),\\ 
\textrm{ with }F:= h^{-1}r^{1-\frac{\delta}{2}}
(ih^{-1}[A,B] -E)r^{1-\frac{\delta}{2}}+\frac{1}{\eps}(A^2+B^2).
\end{gathered}
\end{equation}
We deduce from \eqref{symbol[A,B]} and \eqref{{a,b}} the following 
\begin{lemma}\label{propF}
The operator $F$  is a semiclassical differential operator in the class $S^0(\cjg\xi\cjd^4)$ with semiclassical principal symbol 
\[ \sigma (F)(\xi)= \frac{4|d\varphi|^2}{\eps}+\frac{1}{\eps}(|\xi|^2-|d\varphi|^2)^2+\frac{4}{\eps}(\cjg\xi,d\varphi\cjd)^2.\] 
\end{lemma}
By the semiclassical G{\aa}rding estimate, we obtain the 
\begin{coro}\label{garding}
The operator $F$ of Lemma \ref{propF} is such that there is a constant $C$ so that
\[ \cjg Fu,u\cjd \geq \frac{C}{\eps} (||u||^2_{L^2} +h^2||du||^2_{L^2}).\]
\end{coro}
\textsl{Proof}. It suffices to use that $\sigma(F)(\xi)\geq \frac{C'}{\eps}(1+|\xi|^4)$ for some $C'>0$ and use the 
semiclassical G{\aa}rding estimate. 
\qed

So by writing $\cjg i[A,B]u,u\cjd=\cjg ir^{1-\frac{\delta}{2}}[A,B]r^{1-\frac{\delta}{2}} r^{-1+\frac{\delta}{2}}u,r^{-1+\frac{\delta}{2}}u\cjd$ in \eqref{A+iB} and using \eqref{simpl} and Corollary \ref{garding}, we obtain
that there exists $C>0$ such that for all $u\in C_0^\infty(E_i)$
\begin{equation} \label{phu}
\begin{split}
||P_hu||_{L^2}^2\geq  &\cjg (A^2+B^2)u,u\cjd + \frac{Ch^2}{\eps}(||r^{-1+\frac{\delta}{2}}u||^2_{L^2} +h^2||r^{-1+\frac{\delta}{2}}du||^2_{L^2})+ 
h\cjg Eu,u\cjd\\
& - \frac{h^2}{\eps}(||A (r^{-1+\frac{\delta}{2}} u)||_{L^2}^2 +||B (r^{-1+\frac{\delta}{2}} u)||_{L^2}^2).
\end{split}\end{equation}
We observe that  $h^{-1}[A,r^{-1+\frac{\delta}{2}}]r^{1+\frac{\delta}{2}}\in 
S^0(\cjg \xi\cjd)$ and  $h^{-1}[B,r^{-1+\frac{\delta}{2}}]r^{1+\frac{\delta}{2}}\in 
h S^0(1)$, and thus 
\[ ||A (r^{-1+\frac{\delta}{2}} u)||_{L^2}^2 +||B (r^{-1+\frac{\delta}{2}} u)||_{L^2}^2)\leq C'(||Au||_{L^2}^2+||Bu||_{L^2}^2+
h^2||r^{-1+\frac{\delta}{2}}u||^2_{L^2}+h^4||r^{-1+\frac{\delta}{2}}du||^2_{L^2})\]
for some $C'>0$. Taking $h$ small, this implies with \eqref{phu} that there exists a new constant $C>0$ such that 
\begin{equation}\label{phu2}
||P_hu||_{L^2}^2\geq  \frac{1}{2}\cjg (A^2+B^2)u,u\cjd + \frac{Ch^2}{\eps}(||r^{-1+\frac{\delta}{2}}u||^2_{L^2} +h^2||r^{-1+\frac{\delta}{2}}du||^2_{L^2})+ 
h\cjg Eu,u\cjd.\end{equation}
It remains to deal with $h\cjg Eu,u\cjd$: we first write $E=4|d\varphi_\eps|^2(c(z)(A+h^2\la^2)+
{\rm Op}_h(\ell)B)+hr^{-1+\frac{\delta}{2}}Sr^{-1+\frac{\delta}{2}}$ where $S$
is a semiclassical differential operator in the class $S^0(\cjg\xi\cjd)$ by the decay estimates on $c(z),\ell(z,\xi)$ as $z\to\infty$, 
then by Cauchy-Schwartz (and with $L:={\rm Op}_h(\ell)$)
\[\begin{split} 
|\cjg hEu,u\cjd|\leq & Ch(||Au||_{L^2}+h^2||r^{-1+\frac{\delta}{2}}u||_{L^2}+ h||Sr^{-1+\frac{\delta}{2}}u||_{L^2}) ||r^{-1+\frac{\delta}{2}}u||_{L^2}
+Ch||Bu||_{L^2}||Lu||_{L^2}\\
\leq & \frac{1}{4}||Au||^2_{L^2}+ h^2||Sr^{-1+\frac{\delta}{2}}u||^2_{L^2}+Ch^2||r^{-1+\frac{\delta}{2}}u||_{L^2}^2+\frac{1}{4} ||Bu||^2_{L^2}+ 
Ch^2||Lu||^2_{L^2} 
\end{split}\]
where $C$ is a constant independent of $h,\eps$ but may change from line to line.
Now we observe that $Lr^{1-\frac{\delta}{2}}$ and $S$ are in $S^0(\cjg \xi\cjd)$ and thus 
\[||Sr^{-1+\frac{\delta}{2}}u||^2_{L^2}+ ||Lu||^2_{L^2}\leq C (||r^{-1+\frac{\delta}{2}}u||^2_{L^2}+ h^2||r^{-1+\frac{\delta}{2}}du||^2_{L^2}),\] 
which by \eqref{phu2} implies that there exists $C>0$ such that for all $\eps\gg h>0$ with $\eps$ small enough
\[||P_hu||_{L^2}^2\geq  \frac{Ch^2}{\eps}(||r^{-1+\frac{\delta}{2}}u||^2_{L^2} +h^2||r^{-1+\frac{\delta}{2}}du||^2_{L^2})
\]
for all $u\in C_0^\infty(E_i)$ . The proof is complete.\qed\\

Combining now Proposition \ref{carlemaninend} and \eqref{carlemaninK}, we obtain 
\begin{prop}
\label{carlemanestimate}
Let $(M_0,g_0)$ be a Riemann surface with Euclidean ends with $x$ a boundary defining function of the radial compactification $\bbar{M}_0$ and let $\varphi_\eps=\varphi-\frac{h}{\eps}\varphi_0$ where $\varphi$ is a harmonic function
with non-degenerate critical points and linear growth on $M_0$ 
and $\varphi_0$ satisfies $\Delta_{g_0}\varphi_0=x^{2-\delta}$ as above.
Then for all $V\in x^{1-\frac{\delta}{2}}L^\infty(M_0)$ there exists an $h_0>0$, $\eps_0$ and $C>0$  such that for all $0<h<h_0$, $h\ll \eps<\eps_0$ and $u\in C^\infty_0(M_0)$, we have 
\begin{equation}\label{carleman}
\frac{1}{h}\|x^{1-\frac{\delta}{2}}u\|_{L^2}^2 + \frac{1}{h^2}\|x^{1-\frac{\delta}{2}} u |d\varphi|\|_{L^2}^2 + \|x^{1-\frac{\delta}{2}}du\|_{L^2}^2 \leq C\eps \|e^{\varphi_\eps/h}
(\Delta_g + V-\la^2)e^{-\varphi_\eps/h} u\|_{L^2}^2
\end{equation}
\end{prop}
\noindent\textbf{Proof}.  
As in the proof of Proposition 3.1 in \cite{GT}, by taking $\eps$ small enough, we see that the combination of \eqref{carlemaninK} and Proposition \ref{carlemaninend} shows that for any 
$w \in C_0^\infty(M_0)$, 
\[\begin{gathered}
\frac{C}{\eps}\Big(\frac{1}{h}\|x^{1-\frac{\delta}{2}}w\|_{L^2}^2 + \frac{1}{h^2}\|x^{1-\frac{\delta}{2}}w |d\varphi|\|_{L^2}^2 
+\frac{1}{h^2}\|x^{1-\frac{\delta}{2}}w |d\varphi_{\eps}|\|_{L^2}^2 +
\|x^{1-\frac{\delta}{2}}dw\|_{L^2}^2\Big) \\
\leq \|e^{\frac{\varphi_{\eps}}{h}}(\Delta-\lambda^2)
e^{-\frac{\varphi_\eps}{h}}w\|_{L^2}^2
\end{gathered}
\] 
which ends the proof.\qed

\subsection{The quadratic weight case for surfaces}
In this section, $\varphi$ has quadratic growth at infinity, which corresponds to the case where 
$V\in e^{-\gamma/x^2}L^\infty$ for all $\gamma>0$.
The proof when $\varphi$ has quadratic growth at infinity is even simpler than the linear growth case. We define $\varphi_0\in x^{-2}L^\infty$ 
to be a solution of $\Delta_{g_0}\varphi_0=1$, this is possible by Lemma \ref{rightinv} and one easily obtains 
from Proposition 5.61 of \cite{APS} that $\varphi_0=-x^{-2}/4+O(x^{-1})$ as $x\to 0$. We let 
$\varphi_\eps:=\varphi-\frac{h}{\eps}\varphi_0$ which satisfies $\Delta_{g_0}\varphi_\eps/h=-1/\eps$.

If $K\subset M_0$ is a compact set, the Carleman estimate \eqref{carlemaninK} in $K$
is satisfied by Proposition 3.1 of \cite{GT}, it then remains to get the estimate in the ends $E_1,\dots,E_N$. But the 
exact same proof as in Lemma 3.1 and Lemma 3.2 of \cite{GT} gives directly that for any $w \in C_0^\infty(E_i)$
\begin{equation}\label{CEinEi}
\frac{C}{\eps}\Big(\frac{1}{h}\|w\|_{L^2}^2 + \frac{1}{h^2}\|w |d\varphi|\|_{L^2}^2 +\frac{1}{h^2}\|w |d\varphi_{\eps}|\|_{L^2}^2 +
\|dw\|_{L^2}^2\Big) \leq \|e^{\varphi_{\eps}/h}\Delta_{g_0} e^{-\varphi_{\eps}/h} w\|_{L^2}^2
\end{equation}
for some $C>0$ independent of $\eps,h$ and it suffices to glue the estimates in $K$ and in the ends $E_i$ as in Proposition 3.1 of \cite{GT}, 
to obtain \eqref{CEinEi} for any $w \in C_0^\infty(M_0)$.
Then by using triangle inequality 
\[||e^{\varphi_{\eps}/h}(\Delta_{g_0}+V-\la^2)e^{-\varphi_{\eps}/h}u||_{L^2}\leq ||e^{\varphi_{\eps}/h}\Delta_{g_0}e^{-\varphi_{\eps}/h}u||_{L^2}+
C||u||_{L^2}\]
for some $C$ depending on $\la, ||V||_{L^\infty}$, we see that the $V-\la^2$ term can be absorbed by the 
left hand side of \eqref{CEinEi} and we finally deduce
\begin{prop}\label{CEquad}
Let $(M_0,g_0)$ be a Riemann surface with Euclidean ends  and let $\varphi_\eps=\varphi-\frac{h}{\eps}\varphi_0$ where $\varphi$ is a harmonic function
with non-degenerate critical points and quadratic growth on $M_0$ and 
$\varphi_0$ satisfies $\Delta_{g_0}\varphi_0=1$ with $\varphi_0\in x^{-2}L^\infty(M_0)$.
Then for all $V\in L^\infty$ there exists an $h_0>0$, $\eps_0$ and $C>0$  such that for all $0<h<h_0$, $h\ll \eps<\eps_0$ and $u\in C^\infty_0(M_0)$
\[\frac{C}{\eps}\Big(\frac{1}{h}\|u\|_{L^2}^2 +\frac{1}{h^2}||u |d\varphi| ||^2_{L^2}+
\|du\|_{L^2}^2\Big) \leq \|e^{\varphi_{\eps}/h}(\Delta_{g_0} +V-\la^2)e^{-\varphi_{\eps}/h} u\|_{L^2}^2.\]
\end{prop}
The main difference with the linear weight case is that one can use a convexification which has quadratic growth at infinity which allows to absorb the $\la^2$ term, while it was not the case for the linearly growing weights. 

\end{section}

\begin{section}{Complex Geometric Optics on a Riemann Surface with Euclidean ends}\label{CGOriemann}
As in \cite{Bu,IUY,GT}, the method for identifying the potential at a point $p$ is to construct complex geometric optic solutions depending on a small parameter $h>0$, with phase a Morse holomorphic function with a non-degenerate 
critical point at $p$, and then to apply the stationary phase method. 
Here, in addition, we need the phase to be of linear growth at infinity if $V\in e^{-\gamma/x}L^\infty$ for all $\gamma>0$
while the phase has to be of quadratic growth at infinity if $V\in e^{-\gamma/x^2}L^\infty$ for all $\gamma>0$.\\
 
We shall now assume that $M_0$ is a non-compact surface with genus $g$ with $N$ ends equipped with a metric $g_0$ 
which is Euclidean in the ends, and $V$ is a $C^{1,\alpha}$ function in $M_0$. Moreover, if $V\in e^{-\gamma/x}L^\infty$ for all $\gamma>0$, 
we ask that $N\geq \max(2g+1,2)$ while if $V\in e^{-\gamma/x^2}L^\infty$ for all $\gamma>0$, we assume that $N\geq g+1$. 
As above, let us use a smooth positive function $x$ which is equal to $1$ in a large compact set of $M_0$ and is equal to
$x=|z|^{-1}$ in the regions $|z|>r_0$ of the ends $E_i\simeq \{z\in\cc; |z|>1\}$, where $r_0$ is a fixed large number.
This function is a boundary defining function of the radial 
compactification of $M_0$ in the sense of Melrose \cite{APS}. 
To construct the complex geometric optics solutions, we 
will need to work with the weighted spaces $x^{-\alpha}L^2(M_0)$ where $\alpha\in\rr_+$.\\

Let $\mc{H}$ be the finite dimensional complex vector space defined in the beginning of Section \ref{morseholo}.
Choose $p\in M_0$ such that there exists a Morse holomorphic function $\Phi=\varphi+i\psi\in\mc{H}$ on $M_0$, 
with a critical point at $p$; there is a dense set of such points by Proposition \ref{criticalpoints}. 
The purpose of this section is to construct solutions $u$ on $M_0$ of $(\Delta -\la^2+V)u = 0$ of the form
\begin{equation}
\label{cgo}
u = e^{\Phi/h}(a + r_1 + r_2)
\end{equation}
for $h>0$ small, where $a\in x^{-J+1}L^2$ with $J\in\rr_+\setminus \nn$
 is a holomorphic function on $M_0$, obtained by Lemma \ref{amplitude}, such that $a(p)\not=0$ 
and $a$ vanishing to order 
$L$ (for some fixed large $L$) at all other critical points of $\Phi$, and finally $r_1,r_2$ will be remainder terms which
are small as $h\to 0$ and have particular properties near the critical points of $\Phi$.
More precisely, $e^{\varphi_0/\eps}r_2$ will be a $o_{L^2}(h)$  and $r_1$ will be a $O_{x^{-J}L^2}(h)$ but with 
an explicit expression, which can be used to obtain sufficient information 
in order to apply the stationary phase method.

\subsubsection{Construction of $r_1$}\label{constr1}
We want to construct $r_1=O_{x^{-J}L^2}(h)$ which satisfies
\[e^{-\Phi/h}(\Delta_{g_0} -\la^2+V)e^{\Phi/h}(a + r_1) = O_{x^{-J}L^2}(h) \] 
for some large $J\in\rr_+\setminus \nn$  so that $a\in x^{-J+1}L^2$.

Let $G$ be the operator of Lemma \ref{rightinv}, mapping continuously $x^{-J+1}L^2(M_0)$ 
to $x^{-J-1}L^2(M_0)$. Then clearly $\bar{\partial}\partial G=\frac{i}{2}\star^{-1}$  when acting on $x^{-J+1}L^2$,
here $\star^{-1}$ is the inverse of $\star$ mapping functions to $2$-forms.
First, we will search for $r_1$ satisfying
\begin{equation}
\label{dequation}
e^{-2i\psi/h}\partial e^{2i\psi/h} r_1 = -\pl G (a(V-\la^2)) + \omega + O_{x^{-J}H^1}(h)
\end{equation}
 with $\omega\in x^{-J}L^2(M_0)$ a holomorphic 1-form on $M_0$  
and $\|r_1\|_{x^{-J}L^2} = O(h)$. 
Indeed, using the fact that $\Phi$ is holomorphic we have
\[e^{-\Phi/h}\Delta_{g_0}e^{\Phi /h}=-2i\star  \bar{\pl} e^{-\Phi/h}\pl e^{\Phi/h}=-2i\star  \bar{\pl} e^{-\frac{1}{h}(\Phi-\bar{\Phi})}\pl e^{\frac{1}{h}(\Phi-\bar{\Phi})}=
-2i\star \bar{\pl}e^{-2i\psi/h}\pl e^{2i\psi/h}\]
and applying $-2i\star\bar{\pl}$ to \eqref{dequation}, this gives
\[e^{-\Phi/h}(\Delta_{g_0}+V)e^{\Phi/h}r_1=-a(V-\la^2)+O_{x^{-J}L^2}(h).\]
Writing $-\pl G(a(V-\la^2))=:c(z)dz$ in local complex coordinates, $c(z)$ is $C^{2,\alpha}$ by elliptic regularity and 
we have $2i\pl_{\bar{z}}c(z)=a(V-\la^2)$, therefore $\pl_z\pl_{\bar{z}}c(p')=\pl^2_{\bar{z}}c(p')=0$ at each critical point $p'\not=p$ by construction 
of the function $a$.  Therefore, we deduce that at each critical point $p'\neq p$, $c(z)$ 
has Taylor series expansion $\sum_{j = 0}^2 c_j z^j + O(|z|^{2+\alpha})$. That is, all the lower order terms of the Taylor expansion of $c(z)$ 
around $p'$ are polynomials of $z$ only. By Lemma \ref{control the zero}, and possibly by taking $J$ larger, there exists a holomorphic function $f\in x^{-J}L^2$ 
such that $\omega:=\pl f$ 
has Taylor expansion equal to that of $\pl G(a(V-\la^2))$ at all critical points $p'\not=p$ of $\Phi$. We deduce that, if
 $b:=-\pl G(a(V-\la^2))+\omega=b(z)dz$, we have
\begin{equation}\label{decayofb}
\begin{gathered}
|\pl_{\bar{z}}^m\pl^{\ell}_z b(z)|=O(|z|^{2+\alpha-\ell-m}) , \quad \textrm{ for } \ell+m\leq 2 , \textrm{ at critical points }p'\not=p\\
|b(z)|=O(|z|) , \qquad \qquad \qquad \qquad \textrm{ if }p'=p.  
 \end{gathered}\end{equation}
 Now, we let $\chi_1\in C_0^\infty(M_0)$ be a cutoff function supported in a small neighbourhood $U_p$ of the critical point $p$ and identically $1$ near 
 $p$, and 
$\chi\in C_0^\infty(M_0)$ is defined similarly with $\chi =1$ on the support of $\chi_1$.
We will construct $r_1$ to be a sum $r_1=r_{11} +h r_{12}$ where  $r_{11}$ is a compactly supported 
approximate solution of \eqref{dequation}  near the critical point $p$ of $\Phi$ and $r_{12}$ is correction term supported 
away from $p$.
We define locally in complex coordinates centered at $p$ and containing the support of $\chi$
\begin{equation}\label{defr11}
r_{11}:=\chi e^{-2i\psi/h}R(e^{2i\psi/h}\chi_1b)
\end{equation}
where $Rf(z) := -(2\pi i)^{-1}\int_{\R^2} \frac{1}{\bar{z}-\bar{\xi}}f d\bar{\xi}\wedge d\xi$ for $f\in L^\infty$ compactly supported 
is the classical Cauchy operator inverting locally $\pl_z$ ($r_{11}$ is extended by $0$ outside the neighbourhood of $p$).
The function $r_{11}$ is in $C^{3,\alpha}(M_0)$ and we have 
\begin{equation}\begin{gathered}\label{r11}
e^{-2i\psi/h}\pl(e^{2i\psi/h}r_{11}) = \chi_1(-\pl G(a(V-\la^2)) + \omega) + \eta\\
\textrm{ with }\eta:= e^{-2i\psi/h}R(e^{2i\psi/h}\chi_1b)\pl\chi.
\end{gathered}
\end{equation}
We then construct $r_{12}$ by observing that  $b$ vanishes to order $2+\alpha$ at critical points of $\Phi$ other than $p$ 
(from \eqref{decayofb}), and 
$\pl \chi=0$ in a neighbourhood of any critical point of $\psi$, so we can find $r_{12}$ satisfying
\begin{equation}\label{defr12}
2ir_{12}\pl\psi = (1-\chi_1)b .
\end{equation} 
This is possible since both $\pl\psi$ and the right hand side are valued in $T^*_{1,0}M_0$ and 
$\pl \psi$ has finitely many isolated zeroes on $M_0$: 
$r_{12}$ is then a function which is in $C^{2,\alpha}(M_0\setminus{P})$ where $P:=\{p_1,\dots, p_n\}$ is the set of critical points other than $p$,
it extends to a function in $C^{1,\alpha}(M_0)$  and it satisfies in local complex coordinates $z$ at each $p_j$ 
\[ |\pl_{\bar{z}}^\beta\pl_z^\gamma r_{12}(z)|\leq C|z|^{1+\alpha-\beta-\gamma} , \quad \beta+\gamma\leq 2\]
by using also the fact that $\pl \psi$ can be locally be considered as a smooth function with a zero of order $1$ at each $p_j$.
Moreover $b\in x^{-J}H^2(M_0)$ thus  $r_1\in x^{-J}H^2(M_0)$ and  we have 
\[ e^{-2i\psi/h}\pl(e^{2i\psi/h}r_1) = b+h\pl r_{12}+\eta=-\pl G(a(V-\lambda^2))+\omega+ h\pl r_{12} + \eta.\]

\begin{lemma}\label{fewestimates}
The following estimates hold true 
\[\begin{gathered} 
||\eta||_{H^2(M_0)}=O(|\log h|),\quad \|\eta\|_{H^1(M_0)}\leq O(h|\log h|), \quad ||x^{J}\pl r_{12}||_{H^1(M_0)}=O(1),\\
 ||x^Jr_{1}||_{L^2}=O(h), \quad ||x^J(r_1-h\til{r}_{12})||_{L^2}=o(h) 
\end{gathered} \]
where $\til{r}_{12}$ solves $2i\til{r}_{12}\pl\psi = b$.
\end{lemma}
\textbf{Proof}. The proof is exactly the same as the proof of Lemma 4.2 in \cite{GT2}, except that one needs to add the weight 
$x^J$ to have bounded integrals.
\qed\\

As a direct consequence, we have 
\begin{coro}\label{corerrorterm}
With $r_1=r_{11}+h r_{12}$, there exists $J > 0$ such that 
\[||e^{-\Phi/h}(\Delta_{g_0}-\la^2+V)e^{\Phi/h}(a + r_1)||_{x^{-J}L^2(M_0)} = O(h|\log h|). \]
\end{coro}

\subsubsection{Construction of $r_2$}
In this section, we complete the construction of the complex geometric optic solutions. We deal with the general case of surfaces 
and we shall show the following 
\begin{proposition}
\label{completecgo}
If $\varphi_0$ is the subharmonic function constructed in Section \ref{Carleman}, then for $\eps$ small enough there exist solutions to $(\Delta_{g_0} -\la^2+V)u = 0$ of the form $u=e^{\Phi/h}(a+r_1+r_2)$ 
with $r_1=r_{11}+hr_{12}$  constructed in the previous section and 
$r_2\in e^{-\varphi_0/\eps}L^2$ satisfying $\|e^{\varphi_0/\eps}r_2\|_{L^2}\leq Ch^{3/2}|\log h|$.
\end{proposition}
This is a consequence of the following Lemma (which follows from the Carleman estimate obtained in Section \ref{Carleman} above) 
\begin{lemma}
\label{standardargument}
Let $\delta \in (0,1)$, $V\in x^{1-\frac{\delta}{2}}L^\infty(M_0)$, and $\varphi_\eps=\varphi-\frac{h}{\eps}\varphi_0$ a weight with linear growth at infinity 
as in Proposition \ref{carlemanestimate}.
For all $f\in L^2(M_0)$ and all $h>0$ small enough, there exists a solution $v\in L^2(M_0)$ to the equation
\begin{equation}\label{solvab}
e^{-\varphi_\eps/h}(\Delta_g -\la^2+V) e^{\varphi_\eps/h}v = x^{1-\frac{\delta}{2}}f
\end{equation}
satisfying
\[\|v\|_{L^2(M_0)} \leq Ch^\demi\|f\|_{L^2(M_0)}.\]
If $\varphi_\eps$ has quadratic growth at infinity, the same result is true when $V \in L^{\infty}(M_0)$ but $x^{1-\frac{\delta}{2}}f$ can be replaced by $f\in L^2$ in \eqref{solvab}.
\end{lemma}
\noindent{\bf Proof}. The proof is based on a duality argument. Let $P_h:=e^{\varphi_\eps/h}(\Delta_g -\la^2+V) e^{-\varphi_\eps/h}$
and for all $h >0$ the real vector space $\mc{A}:=\{u\in x^{-1+\frac{\delta}{2}}H^1(M_0); 
P_hu\in L^2(M_0)\}$
equipped with the real scalar product 
\[(u,w)_{\mc{A}}:=\cjg P_hu,P_hw\cjd_{L^2}.\] 
By the Carleman estimate of Proposition \ref{carlemanestimate}, the space $\mc{A}$ is a Hilbert space equipped with the scalar product above if $h < h_0$, and
thus the linear functional $L:w\to \int_{M_0}x^{1-\frac{\delta}{2}}fw \,{\rm dvol}_{g_0}$ on $\mc{A}$ is continuous with norm bounded by 
$Ch^\demi||f||_{L^2}$ by Proposition \ref{carlemanestimate},  
and by Riesz theorem there is an element $u\in\mc{A}$ such that $(.,u)_{\mc{A}}=L$  
and with norm bounded by the norm of 
$L$. It remains to take $v:=P_hu$ which solves $P_h^*v=x^{1-\frac{\delta}{2}}f$  where 
$P_h^*=e^{-\varphi_\eps/h}(\Delta_g -\la^2+V) e^{\varphi_\eps/h}$ is the adjoint of $P_h$ 
and $v$ satisfies the desired norm estimate. The proof when the weight $\varphi_\eps$ has  quadratic growth at infinity  
is the same, but improves
slightly due to the  Carleman estimate of Proposition \ref{CEquad}.
\qed\\

\noindent{\bf Proof of Proposition \ref{completecgo}}. We first solve the equation 
\[(\Delta +V-\la^2)e^{\varphi_\eps/h}\til{r}_2 =  x^{1-\frac{\delta}{2}}e^{\varphi_\eps/h} \Big(x^{-1+\frac{\delta}{2}}e^{-\varphi_\eps/h}(\Delta +V-\la^2)e^{\Phi/h}(a + r_1)\Big)\]
by using Lemma \ref{standardargument} and the fact that for $J$ large, there is $C>0$ such that for all $h<h_0$
\[ ||x^{-1+\frac{\delta}{2}}e^{-\varphi_\eps/h}(\Delta +V-\la^2)e^{\Phi/h}(a + r_1)||_{L^2}\leq C|| x^Je^{-\Phi/h}(\Delta -\la^2+V)e^{\Phi/h}(a + r_1)||_{L^2}\]
since $x^{-J-1}e^{\varphi_0/\eps}\in L^\infty(M_0)$ for all $J$ (recall that $\varphi_0\sim -x^{-\delta}/\delta^2$ as $x\to 0$). 
But now the right hand side is bounded by $O(h|\log h|)$ according to Corollary 
\ref{corerrorterm}, therefore we set $r_2:=-e^{-i\psi/h-\varphi_0/\eps}\til{r}_2$ which 
satisfies $(\Delta_{g_0}-\la^2+V)e^{\Phi/h}(a+r_1+r_2)=0$ and, by Lemma \ref{standardargument}, 
the norm estimate $||e^{\varphi_0/\eps}r_2||_{L^2}\leq O(h^{3/2}|\log h|)$.
 \qed
\end{section}

\begin{section}{Scattering on surface with Euclidean ends} \label{sec_scattering}
Let $(M_0,g_0)$ be a surface with Euclidean ends and  $V\in  e^{-\gamma/x}L^\infty(M_0)$ for some 
$\gamma$. The scattering theory in this setting is described for instance in Melrose \cite{MelStanford}, here 
we will follow this presentation (see also Section 3 in Uhlmann-Vasy \cite{UhlVa} for the $\rr^n$ case).  
First, using standard methods in scattering theory, we define the resolvent on the continuous spectrum as follows
\begin{lemma}\label{resolvent}
The resolvent $R_V(\la):=(\Delta_{g_0}+V-\la^2)^{-1}$ admits a meromorphic extension from $\{{\rm Im}(\la)<0\}$ 
to  $\{{\rm Im}(\la)\leq A, {\rm Re}(\la)\not= 0\}$, as a family of operators mapping $e^{-\gamma/x}L^2(M_0)$ to 
$e^{\gamma/x}L^2(M_0)$ for any $\gamma>A$.  Moreover, for $\la\in\rr\setminus\{0\}$ not a pole,
$R_V(\la)$ maps continuously $x^\alpha L^2$ to $x^{-\alpha}L^2$ for any $\alpha>1/2$.
\end{lemma}
\noindent\textsl{Proof}. The statement is known for $V=0$ and $M_0=\rr^2$ by using the explicit formula of the 
resolvent convolution kernel on $\rr^2$ in terms of Hankel functions (see for instance \cite{MelStanford}),
we shall denote $R_0(\la)$ this continued resolvent. 
More precisely, for all $A>0$, the operator $R_0(\la)$ continues analytically from  $\{{\rm Im}(\la)<0\}$ to 
$\{{\rm Im}(\la)\leq A, {\rm Re}(\la)\not=0\}$ as a family of bounded operators mapping $e^{-\gamma/x}L^2$ to $e^{\gamma/x}L^2$
for any $\gamma>A$.
Now we can set $\chi\in C_0^\infty(M_0)$ such that $1-\chi$ is supported in the ends $E_i$, and let
$\chi_0,\chi_1\in C_0^\infty(M_0)$ such that $(1-\chi_0)=1$ on the support of $(1-\chi)$ and $\chi_1=1$ on the support of $\chi$. 
Let $\la_0\in -i\rr_+$ with $i\la_0\gg 0$, then the resolvent $R_0(\la_0)$ is well defined from $L^2(M_0)$ to $H^2(M_0)$ since the Laplacian is essentially self-adjoint \cite[Proposition 8.2.4]{T2}, and we have a parametrix 
\[E(\la):= (1-\chi_0)R_{0}(\la)(1-\chi)+\chi_1 R_0(\la_0)\chi \] 
which satisfies 
\[\begin{gathered} 
(\Delta_{g_0}-\la^2+V)E(\la)=1+K(\la), \\ K(\la):=([\Delta_{g_0},\chi_1]- (\la^2-\la_0^2)\chi_1)R_0(\la_0)\chi- 
[\Delta_{g_0},\chi_0]R_0(\la)(1-\chi)+VE(\la),
\end{gathered}\]
where here we use the notation $R_0(\la)$ for an integral kernel on $M_0$, which in the charts 
$\{z\in\rr^2; |z|>1\}$ corresponding the ends $E_1,\dots E_N$, is given by the integral kernel of $(\Delta_{\rr^2}-\la^2)^{-1}$. 
 Using the explicit expression of the convolution kernel of $R_0(\la)$ in the ends 
 (see for instance Section 1.5 of \cite{MelStanford}) and the decay assumption
on $V$, it is direct to see that for ${\rm Im}(\la)<A, {\rm Re}(\la)\not=0$, the map $\la \mapsto K(\la)$ a is compact analytic family of bounded operators 
from $e^{-\gamma/x}L^2$ to $e^{-\gamma/x}L^2$ for any $\gamma>A$. Moreover $1+K(\la_0)$ is invertible since $||K(\la_0)||_{L^2\to L^2}\leq 1/2$
if $i\la_0$ is large enough.
Then by analytic Fredholm theory, the resolvent $R_V(\la)$ has an meromorphic extension to ${\rm Im}(\la)<A, {\rm Re}(\la)\not=0$
as a bounded operator from $e^{-\gamma/x}L^2$ to $e^{\gamma/x}L^2$ if $\gamma>A$, given by 
\[ R_V(\la)=E(\la)(1+K(\la))^{-1}.\]
Now $(1+K(\la))^{-1}=1+Q(\la)$ for some $Q(\la)=-K(\la)(1+K(\la))^{-1}$ mapping  $e^{-\gamma/x}L^2$ to itself
for any $\gamma>A$, which proves the mapping properties of $R_V(\la)$ on exponential weighted spaces. For the 
mapping properties on $\{{\rm Re}(\la)=0\}$, a similar argument works.
\qed\\

A corollary of this Lemma is the mapping property
\begin{coro}\label{mapping}
For $ \la\in\rr\setminus \{0\}$ not a pole of $R_V(\la)$, 
and $f\in e^{-\gamma/x}L^\infty$ for some $\gamma>0$, then there exists $v\in C^\infty(\pl\bbar{M}_0)$ such that
\[R_V(\la)f -x^\demi e^{-i\la/x}v \in L^2.\] 
\end{coro}
\noindent\textsl{Proof}. Using the expression $R_V(\la)=E(\la)(1+Q(\la))$ of the proof of Lemma \ref{resolvent}, it suffices
to know the mapping property of $E(\la)$ on $e^{-\gamma/x}L^2$, but since outside a compact set (i.e.~in the ends)
$E(\la)$ is given by the free resolvent on $\rr^2$, this amounts to proving the statement in $\rr^2$, which is 
well-known: for instance, this is proved for $f\in C_0^\infty(\rr^2)$ in Section 1.7 \cite{MelStanford} but the proof extends easily to 
$f\in e^{-\gamma/x}L^\infty(\rr^2)$ since the only used assumption on $f$ for applying a stationary phase argument 
is actually that the  Fourier transform $\hat{f}(z)$ has a holomorphic extension in a complex neighbourhood of $\rr^2$.  
\qed\\

We also have a boundary pairing, the proof of which is exactly the same as 
\cite[Lemma 2.2]{MelStanford} (see also Proposition 3.1 of \cite{UhlVa}).
\begin{lemma}\label{boundarypairing}
For $\la>0$ and $V\in e^{-\gamma/x}L^\infty(M_0)$, if $u_\pm\in x^{-\alpha}L^2(M_0)$ for some $\alpha>1/2$ and 
$(\Delta_{g_0}-\la^2+V)u_\pm\in x^\alpha L^2(M_0)$ with 
\[ u_+-x^\demi e^{i\la/x}f_{++}-x^\demi e^{-i\la/x}f_{+-} \in L^2, \quad u_--x^\demi e^{i\la/x}f_{-+}-x^\demi e^{-i\la/x}f_{--} \in L^2
\] 
for some $f_{\pm\pm}\in C^\infty(\pl \bbar{M}_0)$, then 
\[\cjg u_+, (\Delta_{g_0}+V-\la^2)u_-\cjd-\cjg (\Delta_{g_0}+V-\la^2)u_+, u_-\cjd=2i\la \int_{\pl \bbar{M}_0}(f_{++}\bbar{f_{-+}}-f_{+-}\bbar{f_{--}}) \]
where the volume form on $\pl\bbar{M}_0\simeq\sqcup_{i=1}^N S^1$ is induced by the metric $x^2g|_{T\pl\bbar{M}_0}$. 
\end{lemma}

As a corollary, the same exact arguments as in Sections 2.2 to 2.5 in \cite{MelStanford} show \footnote{In \cite{MelStanford}, 
a unique continuation is used for Schwartz solutions of $(\Delta+V-\la^2)u=0$ when $V$ is a 
compactly supported potential on $\rr^n$  but the same result is also true in our setting, 
this is a consequence of a standard Carleman estimate.}  
\begin{coro}\label{analytic}
The operator $R_V(\la)$ is analytic on  $\la\in \rr\setminus\{0\}$ as a bounded operator from $x^{\alpha}L^2$ to $x^{-\alpha}L^2$
if $\alpha>1/2$.
\end{coro}

In $\rr^2$ there is a Poisson operator $P_0(\la)$ mapping $C^\infty(S^1)$ to $x^{-\alpha}L^2(\rr^2)$ 
for $\alpha>1/2$, which satisfies that for any $f_+\in C^\infty(S^1)$ there exists $f_-\in C^\infty(S^1)$ such that
\[P_0(\la)f_+ -x^\demi e^{i\la/x}f_+ -x^\demi e^{-i\la/x}f_-\in L^2 , \quad (\Delta-\la^2)P_0(\la)f_+=0.\]
We can therefore define in our case a similar Poisson operator $P_V(\la)$ mapping $C^\infty(\pl\bbar{M}_0)$ to $x^{-\alpha} L^2$ for $\alpha>1/2$, 
 by 
\begin{equation}\label{poisson}
P_V(\la)f_+:= (1-\chi)P_0(\la)f_+- R_V(\la)(\Delta_{g_0}+V-\la^2)(1-\chi)P_0(\la)f_+
\end{equation}
where $1-\chi\in C^\infty(M_0)$ equals $1$ in the ends $E_i$ and $P_0(\la)$ denotes here the Schwartz kernel  of the Poisson operator on $\rr^2$ 
pulled back to each of the Euclidean ends $E_i$ of $M_0$ in the obvious way. 
Then, since $(\Delta_{g_0}+V-\la^2)(1-\chi)P_0(\la)f_+\in e^{-\gamma/x}L^2$ for all $\gamma>0$, it suffices to use 
Corollaries \ref{mapping} and \ref{analytic} to see that it defines an analytic Poisson operator $P_V(\la)$ 
on $\la\in \rr\setminus\{0\}$ satisfying that for all $f_+\in C^\infty(\pl\bbar{M_0})$, there exists $f_-\in C^\infty(\pl\bbar{M}_0)$ such that 
\begin{equation}\label{poissonV}
 P_V(\la)f_+ -x^\demi e^{i\la/x}f_+ -x^\demi e^{-i\la/x}f_-\in L^2 , \quad (\Delta+V-\la^2)P_V(\la)f_+=0.
\end{equation}
Moreover,  it is easily seen to be the unique solution of \eqref{poissonV}: indeed, if two such solutions exist then the difference is 
a solution $u$ with asymptotic $x^{\demi}e^{-i\la/x}f_-+L^2$ for some $f_-\in C^\infty(\pl\bbar{M}_0)$, but applying Lemma \ref{boundarypairing}
with $u_-=u_+=u$ shows that $f_-=0$, thus $u\in L^2$, which implies $u=0$ by Corollary \ref{analytic}.

\begin{definition} 
The scattering matrix $S_V(\la):C^\infty(\pl\bbar{M}_0)\to C^\infty(\pl\bbar{M}_0)$ 
for $\la\in \rr\setminus\{0\}$ is defined to be the map $S_V(\la)f_+:=f_-$ where $f_-$ is given by the asymptotic
\[P_V(\la)f_+= x^\demi e^{i\la/x}f_+ +x^\demi e^{-i\la/x}f_- +g ,\,\, \textrm{ with }\,\, g\in L^2.\]
\end{definition}
We remark that, using Lemma \ref{boundarypairing} and the uniqueness of the Poisson operator, one easily deduces for $\la\in\rr\setminus\{0\}$
\begin{equation}\label{relationSV}
S_V(\la)^*=S_V(-\la)=S_V(\la)^{-1}
\end{equation} 
where the scalar product on $L^2(\pl\bbar{M}_0)$ is induced by the metric $x^2g_0|_{T\pl\bbar{M}_0}$.

We can now state a density result similar to Proposition 3.3 of \cite{UhlVa}: 
\begin{prop}\label{density}
If $V\in e^{-\gamma_0/x}L^\infty(M_0)$ $($resp. $V\in e^{-\gamma_0/x^2}L^\infty(M_0))$ 
for some $\gamma_0>0$, and $\la\in\rr\setminus\{0\}$, then for any $0<\gamma<\gamma'<\gamma_0$ the set  
\[\mc{F}:=\{ P_V(\la)f_+; f_+\in C^\infty(\pl\bbar{M}_0) \}\]
is dense in the null space of $\Delta_{g_0}+V-\la^2$ in $e^{\gamma/x}L^2(M_0)$ for the topology of $e^{\gamma'/x}L^2(M_0)$ 
$($resp.~in $e^{\gamma/x^2}L^2(M_0)$ for the topology of $e^{\gamma'/x^2}L^2(M_0))$.
\end{prop}
\noindent\textsl{Proof}. First assume $V\in e^{-\gamma_0/x}L^\infty(M_0)$. Let $w\in e^{-\gamma'/x}L^2$ be orthogonal to $\mc{F}$, and set $u_-:=R_V(\la)w$  and $u_+=P_V(\la)f_{++}$ for some
$f_{++} \in C^\infty(\pl\bbar{M}_0)$. Then, define $f_{--}\in C^\infty(\pl\bbar{M}_0)$ by $R_V(\la)w-x^\demi e^{-i\la/x}f_{--}\in L^2$,
and from Lemma \ref{boundarypairing} we obtain $\cjg f_{+-},f_{--}\cjd=0$ since $\cjg w,P_V(\la)f_{++}\cjd=0$ by assumption.
Since $f_{+-} = S_V(\lambda) f_{++}$ is arbitrary, then $f_{--}=0$ and $u_-\in L^2$. In particular, from the parametrix constructed in the proof of 
Lemma \ref{resolvent}
\[R_V(\la)w-  (1-\chi_0)R_0(\la)(1-\chi)(1+Q(\la))w  \in L^2\]  
with $(1+Q(\la))w\in e^{-\gamma'/x}L^2$. Since in each end, $R_0(\la)$ is the integral 
kernel of the free resolvent of the Euclidean Laplacian on $\rr^2$ and $(1-\chi_0)$ and $(1-\chi)$ are supported in the ends, 
we can view the term $(1-\chi_0)R_0(\la)(1-\chi)(1+Q(\la))w$ as a disjoint sum (over the ends) 
of functions on $\rr^2$ of the form 
\begin{equation}\label{fouriertr} 
(1-\chi_0(z)) \frac{1}{(2\pi)^{2}} \int_{\rr^2}e^{iz\xi} (\xi^2-\la^2-i0)^{-1}\hat{f}(\xi)d\xi 
\end{equation}
where in each end $E_i$, $f=(1-\chi)(1+Q(\la))w\in e^{-\gamma'/x}L^2(E_i)$ can be considered as a function in $e^{-\gamma'|z|}L^2(\rr^2)$.
By the Paley-Wiener theorem, $\hat{f}$ is holomorphic in a strip $U=\{|{\rm Im}(\xi)|<\gamma'\}$ with  
bound $\sup_{\eta\leq \gamma}|| \hat{f}(\cdot+i\eta)||_{L^2(\rr^2)}<\infty$ for all $\gamma<\gamma'$, so the fact that \eqref{fouriertr} is in $L^2$ implies that 
$\hat{f}$ vanishes at the real sphere $\{\xi\in\rr^2 ; \xi^2=\la^2\}$, and thus there exists $h$ holomorphic in $U$ such that
$\hat{f}(\xi)=(\xi^2-\la^2)h(\xi)$ (see e.g.~the proof of Lemma 2.5 in \cite{PSU}), and satisfying the same types of $L^2$ estimates
as $\hat{f}$ in $U$ on lines ${\rm Im}(\xi)=\textrm{cst}$. By the Paley-Wiener theorem again, we deduce that  
\eqref{fouriertr} is in $e^{-\gamma|z|}L^2$ and thus $R_V(\la)w\in e^{-\gamma/x}L^2(M_0)$ for any $\gamma<\gamma'$.  
Then if $v\in e^{\gamma/x}L^2(M_0)$ and $(\Delta_{g_0}+V-\la^2)v=0$, one has by integration by parts 
\[ 0=\cjg R_V(\la) w,(\Delta_{g_0}+V-\la^2)v\cjd = \cjg w,v\cjd\]
which ends the proof in the case $V\in e^{-\gamma_0/x}L^\infty(M_0)$. The quadratic decay case $V\in e^{-\gamma_0/x^2}L^\infty(M_0)$
is exactly similar but instead of Paley-Wiener theorem, we use Corollary \ref{corapp} and the inclusions 
$e^{-\gamma'/x^2}L^2\subset e^{-\gamma''/x^2}L^1\cap e^{-\gamma''/x^2}L^2$  and $e^{-\gamma'/x^2}L^\infty\subset 
e^{-\gamma/x^2}L^2$ for  all $\gamma<\gamma''<\gamma'$.
\qed

\end{section}

\begin{section}{Identifying the potential} \label{sec_identify}
\subsection{The case of a surface}
On a Riemann surface $(M_0,g_0)$ with $N$ Euclidean ends and genus $g$, we assume that $V_1,V_2\in C^{1,\alpha}(M_0)$ are two real valued potentials such that the respective scattering operators 
$S_{V_1}(\la)$ and $S_{V_2}(\la)$ agree for a fixed $\la>0$.  We also assume that for all $\gamma>0$  
\[V_1,V_2\in \left\{\begin{array}{ll}
e^{-\gamma/x}L^\infty(M_0) & \textrm{ if }N\geq \max(2g+1,2)\\
e^{-\gamma/x^2}L^\infty(M_0) & \textrm{ if }N\geq g+1.
\end{array}\right.\]
By considering the asymptotics of $u_1:=P_{V_1}(\la)f_1$ and $P_{V_2}(-\la)f_2$ for $f_i\in C^\infty(\pl\bbar{M}_0)$
 we easily have by integration by parts that 
\begin{equation}\label{scatident}
\begin{split}
\int_{M_0} (V_1-V_2)u_1\overline{u_2}\, {\rm dvol}_{g_0}=& -2i\la \int_{\pl\bbar{M}_0} S_{V_1}(\la)f_1.\overline{f_2} - f_1.\overline{S_{V_2}(-\la)f_2}\\  
=& -2i\la \int_{\pl\bbar{M}_0} (S_{V_1}(\la)-S_{V_2}(\la))f_1.\overline{f_2}=0
\end{split}\end{equation} 
by using \eqref{relationSV}. From Proposition \ref{density}, this implies by density that, if $V\in e^{-\gamma/x}L^\infty$ 
(resp. $V\in e^{-\gamma/x^2}L^\infty$ for all $\gamma>0$), then for all solutions $u_i$
of $(\Delta_{g_0}+V_i-\la^2)u_i=0$ in $e^{\gamma'/x}L^2(M_0)$ (resp. $u_i\in e^{\gamma'/x^2}L^2(M_0)$) 
for some $\gamma'>0$, we have 
\begin{equation}\label{integralform} 
\int_{M_0} (V_1-V_2)u_1 \overline{u_2} \, {\rm dvol}_{g_0}=0.
\end{equation} 
We shall now use our complex geometric optics solutions as special solutions in the weighted space 
$e^{-\gamma'/hx}L^2(M_0)$ (resp. $e^{-\gamma'/hx^2}L^2(M_0)$)
for some $\gamma'>0$ if $V\in e^{-\gamma/x}L^\infty$ (resp. $V\in e^{-\gamma/x^2}L^\infty$) for all $\gamma>0$.\\

Let $p\in M_0$ be such that, using Proposition \ref{criticalpoints}, we can choose a holomorphic Morse function $\Phi=\varphi+i\psi$ 
with linear or quadratic growth on $M_0$ (depending on the topological assumption), 
with a critical point at $p$. Then for the complex geometric optics solutions $u_1,u_2$ with phase $\Phi$ 
constructed in Section \ref{CGOriemann}, the identity \eqref{integralform} holds true. 
We will then deduce the 
\begin{proposition}
\label{identcritpts}
Let $\la\in (0,\infty)$ and assume that $S_{V_1}(\la) = S_{V_2}(\la)$, then $V_1(p)= V_2(p)$.
\end{proposition}
\noindent{\bf Proof}. 
Let $u_1$ and $u_2$ be solutions on $M_0$ to 
\[(\Delta_g +V_j-\la^2)u_j = 0\]
constructed in Section \ref{CGOriemann} with phase $\Phi$ for $u_1$ and $-\Phi$ for $u_2$, thus of the form
\[u_1 = e^{\Phi/h}(a + r_1^{1} + r_2^{1}), \quad u_2=e^{-\Phi/h}(a+r_1^{2}+r_2^{2}).\]
We have the identity 
\[\int_{M_0}u_1(V_1 - V_2) \bbar{u_2} \,{\rm dvol}_{g_0}=0\] 
Then by using the estimates in Lemma \ref{fewestimates} and Proposition \ref{completecgo} we have, as $h\to 0$,
\[\int_{M_0}e^{2i\psi/h}|a|^2(V_1 - V_2) \,{\rm dvol}_{g_0} +  
h \int_{M_0}e^{2i\psi/h}(\bbar{a} \til{r}_{12}^{1} + a \bbar{\til{r}_{12}^{2}})(V_1 - V_2) \,{\rm dvol}_{g_0} + o(h)  = 0 \]
where $\til{r}_{12}^{j}\in L^\infty(M_0)$ are defined in Lemma \ref{fewestimates}, with the superscript $j$ refering to the solution for the 
potential $V_j$;  in particular these functions $\til{r}^j_{12}$ are independent of $h$.

By splitting $V_i(\cdot)=(V_i(\cdot)-V_i(p))+V_i(p)$ and using the $C^{1,\alpha}$ regularity assumption on $V_i$,
one can use stationary phase for the $V_i(p)$ term and integration by parts to gain a power of $h$ for the $V_i(\cdot)-V_i(p)$ term
(see the proof of Lemma 5.4 in \cite{GT2} for details) to deduce
\[\int_{M_0}e^{2i\psi/h}|a|^2(V_1 - V_2) \,{\rm dvol}_{g_0} = Ch(V_1(p) - V_2(p)) + o(h)\]
for some $C\neq 0$.  Therefore,
\[ Ch(V_1(p) - V_2(p))+  h\int_{M_0}e^{2i\psi/h}(\bbar{a} \til{r}_{12}^{1} + a \bbar{\til{r}_{12}^{2}})(V_1 - V_2) \,{\rm dvol}_{g_0} = o(h).\]
Now to deal with the middle terms, it suffices to apply a Riemann-Lebesgue type argument like Lemma 5.3 of \cite{GT2} 
to deduce that it is a $o(h)$. 
The argument is simply to approximate the amplitude in the $L^1(M_0)$ norm by a smooth compactly supported function 
and then use stationary phase to deal with the smooth function. We have thus proved that $V_1(p)=V_2(p)$ by taking $h\to 0$.  
\qed

\end{section}

\begin{section}{Appendix}

To obtain mapping properties of the resolvent of $\Delta_{\rr^2}$ acting on functions with Gaussian decay, 
we shall give two Lemmas on Fourier transforms of functions with Gaussian decay. 
\begin{lemma}\label{app1}
Let $f(z)\in e^{-\gamma|z|^2}L^2(\rr^2)$ for some $\gamma>0$. Then the Fourier transform $\hat{f}(\xi)$ extends analytically to 
$\cc^2$ and for all $\xi, \eta\in\rr^2$,
\[ ||\hat{f}(\xi+i\eta) ||_{L^2(\rr^2,d\xi)}\leq 2\pi e^{\frac{|\eta|^2}{4\gamma}}||e^{\gamma|z|^2}f||_{L^2(\rr^2)}.\]
If $f(z)\in e^{-\gamma|z|^2}L^1(\rr^2)$ for some $\gamma>0$ then
\[\sup_{\xi\in\rr^2}|\hat{f}(\xi+i\eta) |\leq e^{\frac{|\eta|^2}{4\gamma}}||e^{\gamma|z|^2}f||_{L^1(\rr^2)}.\]
\end{lemma}
\noindent\textsl{Proof}. The first statement is clear. For the bound, we write 
\[ \hat{f}(\xi+i\eta)=e^{\frac{|\eta|^2}{4\gamma}}\int_{\rr^2}e^{-i\xi. z} e^{-\gamma|z-\frac{\eta}{2\gamma}|^2}e^{\gamma|z|^2}f(z)dz=
e^{\frac{|\eta|^2}{4\gamma}}\mc{F}_{z\to\xi} (e^{-\gamma|z-\frac{\eta}{2\gamma}|^2}e^{\gamma|z|^2}f(z)).
\] 
But the function $e^{-\gamma|z-\frac{\eta}{2\gamma}|^2}e^{\gamma|z|^2}f(z)$ is in $L^2(\rr^2,dz)$ and its norm is bounded uniformly
by $||e^{\gamma|z|^2}f||_{L^2}$, thus it suffices to use the Plancherel theorem to obtain the desired bound. The $L^\infty$ bound is similar.
\qed

\begin{lemma}\label{app2}
Let $F(\xi+i\eta)$ be a complex analytic function on $\rr^2+i\rr^2=\cc^2$ such that there is $C>0$ and $\gamma>0$ with
\[||F(\xi+i\eta) ||_{L^2(\rr^2,d\xi)}\leq Ce^{\frac{|\eta|^2}{4\gamma}} \textrm{ and } \sup_{\xi\in\rr^2}|F(\xi+i\eta) |\leq Ce^{\frac{|\eta|^2}{4\gamma}}.\]  
If $F$ vanishes on the real submanifold $\{|\xi|^2=\la^2\}$, then $\mc{F}^{-1}_{\xi\to z}(\frac{F(\xi)}{|\xi|^2-\la^2})\in e^{-\gamma|z|^2}L^\infty(\rr^2)$.
\end{lemma}
\noindent\textsl{Proof}. First by analyticity of $F$, one has that $F$ vanishes on the complex hypersurface $M_\la:=\{\zeta\in\cc^2; \zeta.\zeta=\la^2\}$ 
(see for instance the proof of Lemma 2.5 of \cite{PSU}), and in particular $G(\zeta)=F(\zeta)/(\zeta.\zeta-\la^2)$ is an analytic function on $\cc^2$.
We will first prove that for each $\eta\in\rr^2$, $G(\xi+i\eta)\in L^1(\rr^2,d\xi)\cap L^\infty(\rr^2,d\xi)$ and 
\begin{equation}\label{Gxi}
||G(\xi+i\eta)||_{L^1(\rr^2,d\xi)}\leq Ce^{\frac{|\eta|^2}{4\gamma}}. 
\end{equation}

If $|\eta| \leq 2$ we choose the disc $B := \{ \xi \in \rr^2; |\xi|^2 < 2(4+\la^2) \}$ and let $\zeta:=\xi+i\eta$. Then $||G(\xi+i\eta)||_{L^1(B,d\xi)}$ and $||(\zeta.\zeta-\la^2)^{-1}||_{L^2(\rr^2 \setminus B,d\xi)}$ are uniformly bounded for $|\eta| \leq 2$, and we obtain by Cauchy-Schwarz that \eqref{Gxi} holds for $|\eta| \leq 2$. For the case $|\eta| > 2$ we define $U_\eta:=\{\xi \in \rr^2; |\zeta.\zeta-\la^2|>|\eta|\}$ and note that 
\begin{equation*}
\begin{gathered}
\sup_{|\eta| > 2} ||(\zeta.\zeta-\lambda^2)^{-1}||_{L^1(\rr^2 \setminus U_{\eta},d\xi)} < \infty, \\
\sup_{|\eta| > 2} ||(\zeta.\zeta-\lambda^2)^{-1}||_{L^2(U_{\eta},d\xi)} < \infty.
\end{gathered}
\end{equation*}
These results follow by decomposing the integration sets to parts where one can change coordinates $\xi_1 + i\xi_2$ to $\tilde{\xi}_1 + i \tilde{\xi}_2 := \zeta.\zeta-\lambda^2$, and by evaluating simple integrals. Then \eqref{Gxi} follows from Cauchy-Schwarz and the estimates for $F$.

Let $\eta=2\gamma z$, we use a contour deformation from $\rr^2$ to $2i\gamma z+\rr^2$ in $\cc^2$,
\[\int_{\rr^2}e^{iz.\xi}G(\xi)d\xi= \int_{\rr^2}e^{iz.(\xi+2i\gamma z)}G(\xi+2i\gamma z)d\xi, \]
which is justified by the fact that  $G(\xi+i\eta)\in L^1(\rr^2 \times K,d\xi \,d\eta)$ for any compact set $K$ in $\rr^2$ by the uniform bound \eqref{Gxi}. Now using \eqref{Gxi} again shows that 
\[\Big|\int_{\rr^2}e^{iz.\xi}G(\xi)d\xi\Big|\leq C e^{-\gamma |z|^2}\]
which ends the proof. 
\qed

\begin{coro} \label{corapp}
Let $f(z)\in e^{-\gamma|z|^2}L^2(\rr^2)\cap e^{-\gamma|z|^2}L^1(\rr^2)$ for some $\gamma>0$. Assume that 
its Fourier transform $\hat{f}(\xi)$ vanishes on the sphere $\{|\xi|=|\la|\}$, then one has  
\[ \mc{F}^{-1}_{\xi\to z}\Big( \frac{\hat{f}(\xi)}{|\xi|^2-\la^2}\Big) \in e^{-\gamma|z|^2}L^\infty(\rr^2).\]
\end{coro}

\end{section}

   

\end{document}